 \documentclass[a4paper,11pt]{article}
  \usepackage{amssymb,amsbsy,amsmath,latexsym,theorem,eepic,enumerate}
  \usepackage{xypic,times}
  \xyoption{curve}
  \usepackage{geometry,color}
  \newcommand{\quads}[4]{\left( #1 \hspace{1em}
\overset{\textstyle{#2}}{\underset{\textstyle{#4}}
{\rule{0mm}{1mm}}} \hspace{1em} #3 \right)}

\newcommand{\labto}[1]{\stackrel{#1}{\longrightarrow}}

\def\midsq{\ar @{} [dr] |}

\def\bu{\bullet}

\newcommand{\threeaxes}[3]{\def\objectstyle{\scriptstyle}  \objectmargin={0pt}
\xy
(0,0)*+{}="a",(0,-6)*+{\rule{0em}{1.5ex}#2}="b",(7,0)*+{\;#1}="c",
(14,-3)*+{\;#3}="d" \ar@{->} "a";"b" \ar @{->}"a";"c"  \ar
@{->}"a";"d"\endxy }

\newcommand{\directs}[2]{\def\objectstyle{\scriptstyle}  \objectmargin={0pt}
\xy
(0,4)*+{}="a",(0,-2)*+{\rule{0em}{1.5ex}#2}="b",(7,4)*+{\;#1}="c"
\ar@{->} "a";"b" \ar @{->}"a";"c" \endxy }

\newcommand{\xdirects}[2]{\def\objectstyle{\scriptstyle}  \objectmargin={0pt}
\xy
(0,0)*+{}="a",(0,-6)*+{\rule{0em}{1.5ex}#2}="b",(7,0)*+{\;#1}="c"
\ar@{->} "a";"b" \ar @{->}"a";"c" \endxy }
\newcommand{\sdirects}[2]{\def\objectstyle{\scriptstyle}  \objectmargin={0pt}
\xy
(0,2.2)*+{}="a",(0,-2.5)*+{\rule{0em}{1.5ex}#2}="b",(7,2.2)*+{\;#1}="c"
\ar@{->} "a";"b" \ar @{->}"a";"c" \endxy }

\newcommand{\br}{\mbox{\rule{0.7em}{0.2ex}\hspace{-0.04em}\rule{0.08em}{1.7ex}}}

\newcommand{\tl}{\mbox{\rule{0.08em}{1.7ex}\rule[1.54ex]{0.7em}{0.2ex}}}

\newcommand{\hh}{\mbox{\rule{0.7em}{0.2ex}\hspace{-0.7em}\rule[1.5ex]{0.70em}{0.2ex}}}

\newcommand{\vv}{\mbox{\rule{0.08em}{1.7ex}\hspace{0.6em}\rule{0.08em}{1.7ex}}}

\newcommand{\sq}{\mbox{\rule{0.08em}{1.7ex}\hspace{-0.00em}\rule{0.7em}{0.2ex}\hspace{-0.7em}\rule[1.54ex]{0.7em}{0.2ex}\hspace{-0.03em}\rule{0.08em}{1.7ex}}}

\newcommand{\tsq}{\mbox{\rule{0.04em}{1.55ex}\hspace{-0.00em}\rule{0.7em}{0.1ex}\hspace{-0.7em}\rule[1.5ex]{0.7em}{0.1ex}\hspace{-0.03em}\rule{0.04em}{1.55ex}}}

\def\rho{\varrho}

\def\eps{\varepsilon}
\def\epsilon{\varepsilon}

\def\prt{\partial}

\def\le{\leqslant}

\def\eps{\varepsilon}
\def\epsilon{\varepsilon}

\def\prt{\partial}

\def\le{\leqslant}

 \def\eps{\varepsilon}

\def\epsilon{\varepsilon}

\def\C{\mathsf{C}}

\begin{document}
\title{Category Theory and Higher Dimensional Algebra: \\potential descriptive tools \\
in neuroscience\thanks{This paper is an extended account of a
presentation given  at the International Conference on Theoretical
Neurobiology, Delhi, Feb 24-26, 2003, by the first author, who would
like to thank the National Institute for Brain Research  for
support. Both authors would like to thank Posina Rayadu for
correspondence. The paper appeared in Proceedings of the
International Conference on Theoretical Neurobiology, Delhi,
February 2003, edited by Nandini Singh, National Brain Research
Centre, Conference Proceedings 1 (2003) 80-92. The references have
been updated 09/02/08. }}
\author{Ronald Brown \and Timothy Porter\thanks{Mathematics Division, School of
Informatics, University of Wales, Bangor, Gwynedd LL57 1UT,
U.K..\newline  email: \{r.brown,t.porter\}@bangor.ac.uk
\hspace{2em}
 http://www.bangor.ac.uk/$\sim$\{mas010,mas013\}}}
\date{}
\maketitle

\begin{abstract}
We explain the notion of colimit in category theory as a potential
tool for describing structures and their communication, and the
notion of higher dimensional algebra as potential yoga for dealing
with processes and processes of processes.
\end{abstract}

\section*{Introduction}
There does seem to be a problem in neuroscience in finding a
language suitable for describing brain activity in a way which
could lead to deduction, evaluation of theories, and even perhaps
calculation. This is especially so when dealing with questions of
global activity, as against activity of individual organs or
cells. How can we bridge the gap between neuronal activity and
what are variously described in the literature as percepts,
concepts, thoughts, emotions, ideas, and so on? What are the
relations between the meanings of these various words? When should
we use one rather than another?

One hope is that a mathematics will arise which could help in
these problems. To encourage such a mathematics, we need a
dialogue between mathematicians and neuroscientists.
Mathematicians can contribute by showing the way mathematics
works, describing processes such as abstraction, concept
refinement, etc., explaining  what is currently available, and
analysing  the deficiences of current mathematics in helping with
these problems. It is likely to be a difficult process to move
towards such new mathematics, since life has evolved for a long
period, whereas language and mathematics are relatively recent.
Yet the ability to do mathematics is itself a result of evolution,
and mathematics has a good track record in scientific discovery.

We cannot expect too much very quickly, since the mathematics so
far developed may not be appropriate for this particular task. On
the other hand it is potentially useful to see how mathematics
deals with particular problems, such as the relations between
words, and also to see what is now available which might have
potential uses. Just as the physical sciences of physics and
chemistry  have strongly stimulated the development of areas of
mathematics, so we can hope that the same will eventually hold for
the neurological sciences.

\section{Relations between words}
The Greeks devised the axiomatic method, but thought of it in a
different manner to that we do today. One can imagine that the way
Euclid's Geometry evolved was simply through the delivering of a
course covering the established facts of the time. In delivering
such a course, it is natural to formalise the starting points, and
so arranging a sensible structure. These starting points came to
be called postulates, definitions  and axioms, and they were
thought to deal with real, or even ideal, objects, named points,
lines, distance and so on. The modern view, initiated by the
discovery of non Euclidean geometry,  is that the words points,
lines, etc. should be taken as undefined terms, and that axioms
give the relations between these.  This allows the axioms to apply
to many other instances, and has  led to the power of modern
geometry and algebra.

This suggests a task for the professionals in neuroscience, in
order to help a mathematician struggling with the literature,
namely to devise some kind of glossary with clear relations
between these various words and their usages, in order to see what
kind of axiomatic system is needed to describe their
relationships. Clarifying, for instance, the meaning to be
ascribed to `concept', `percept', `thought', `emotion', etc., and
above all the relations between these words, is clearly a
fundamental but difficult step.

\section{Category theory and colimits: gluing and structure}
One of the strong developments in mathematics of the 20th century
has been that of category theory, with its power of describing the
processes of mathematics, developing new logics, unifying
different topics, and revealing underlying abstract processes
which have turned out to have wide implications and uses.

Abstraction allows analogies,  by encoding relations, and
relations among relations. As an example, when we note that $2 +3=
3+2$ and  $2\times 3 = 3 \times 2$, and extend this to the
abstract  {\it commutative law} $x\circ y =y \circ x$ for a binary
operation $\circ$,  we are making an analogy between addition and
multiplication, and also make this law available in  other
situations. One may presume that the power of abstraction, in some
sense of making maps, must be deeply encoded in evolutionary
history as a technique for encouraging survival, since a map gives
a small and manipulable model of the environment. Symbolic
manipulation in mathematics often involves {\it rewriting}, an
example of which is using the commutative law, i.e. replacing in a
complicated formulae various instances of $x \circ y$ by $y \circ
x$. We do this rewriting for example in obtaining the equation
$$ 3 \times 2 \times 5 \times 3 \times 2 \times 5 \times 3 = 2^2 \times 3 ^3 \times 5^2.$$
For more on rewriting, see \cite{KB}.

A study of why and  how mathematics works could be useful for
making models  for neurological functions involving maps of the
environment. Mathematics may also possibly provide  a
comprehensible case study of the evolution  of complex interacting
structures, and so may yield analogies helpful for developing and
evaluating  models of brain activity, in order to derive better
models, and so better understanding. We expect to need a new
language, a new mathematics,  for describing brain activity. To
see what is involved in this search, it is reasonable to study the
evolution of mathematics, and of particular branches such as
category theory.

Some description of category theory is given in \cite{B-P-cyt}.
The relevance to biological development is described in a series
of papers by Ehresmann and Vanbremersch \cite{E&V}. In particular,
they see the notion of {\it colimit} in a category as describing a
structure made up of inter-related parts, so that a category
evolving with time can then allow evolving structures, the
structures being given as colimits in the category $\C_t$ at time
$t$.

 This notion of colimit gives a very general setting
in which to describe the process of {\it gluing} or {\it
amalgamating} complex structures, together with a description of
the method of input to and output from these structures.

First we need to give the notion of category. This developed from
a useful notation for a function: moving from the somewhat obscure
`a function is $y=f(x)$ where $y$ varies as $x$ varies' to the
clearer `a function  $f: X \to Y$ assigns to each element $x$ of
the set $X$ an element $f(x)$ of the set  $Y$'. This sees
`function' as being a `process'. The composition of functions then
suggests the first step in the notion of a {\bf category}
$\mathsf{C}$, which consists of a class $Ob(\C)$ of `objects', a
set of `arrows', or `morphisms' $f: X \to Y$ for any two objects
$X,Y$, and a composition, giving, for instance,  $fg: X \to Z$ if
also $g: Y \to Z$. This composite is represented by the diagram
\begin{equation}
X \labto{f} Y \labto{g} Z
\end{equation}
or even by
\begin{equation}\xymatrix{
X \ar [r]^f \ar [dr] _{fg}&  Y \ar [d]^g \\
& Z}\label{commdiag}
\end{equation}
A category thus has not only a composition structure but also a
`position' structure given by its class of objects. The only rules
are associativity $f(gh)=(fg)h$ when both sides are defined, and
the existence of identities $1_X$ at each object $X$, so that with
$f$ as above, $1_X f=f = f1_Y$.

The notion of colimit in a category generalises the notion of
forming the union $X\cup Y$ of two overlapping sets,  with
intersection $X \cap Y$. However, instead of concentrating on the
sets $X,Y$ themselves, we place them in context, and say that the
utility of the union  is that it allows us to construct {\it
functions} $f: X \cup Y \to C$ for any $C$ by specifying functions
$f_X: X \to C, f_Y : Y \to C$ which agree on the intersection $W=
X \cap Y$. So we replace the specific construction of $X \cup Y$
by a property which describes, using functions, the relation of
this construction to all other sets. That is, the emphasis is on
the relation between input and output.

A colimit has `input data', a `cocone'. In the case of $X \cup Y$,
this cocone consist of the two functions $i_X: X \cap Y \to X,
\;i_X: X \cap Y \to Y$.

 There are  similar situations in other contexts. The numbers
$$ \max(a,b), \mbox{ and } \mathrm{lcm}(a,b)$$
are all constructed from  their `parts', the numbers $a,b$. Here
the `arrows' of our previous `sets and functions' example  are
replaced in the first case by the order relation `$x \le y$', `$x$
is less than or equal to $y$', and in the second case by the
divisibility relation `$x | y$', `$x$ divides $y$'. So we have the
rules that: (i) if $a \le c$ and $b \le c$ then $\max(a,b)\le c$,
(ii) if $a | c$ and $b | c$ then $\mathrm{lcm}(a,b)| c$.  Thus to
make analogies between constructions for many different
mathematical structures, we simply formulate a notion in a general
category -- that is all! In this way  category theory has been a
great unifying force in mathematics of the 20th century, and
continues so to do.

We also generalise to more complex input data. So the `input data'
for a colimit is   a {\it diagram} $D$, that is a collection of
some objects in a category $\C$ and some arrows between them, such
as:
$$ 
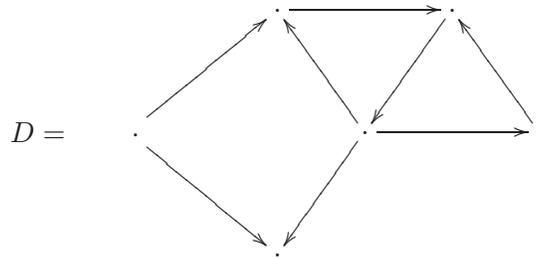
\xymatrix@R=3pc{ & . \ar[rr]&&.\ar[dl]&\\
**[l] D = \qquad \cdot \ar[ur]\ar[dr]& & . \ar[ul]\ar[rr]\ar[dl] & &.\ar[ul]\\
&.&&&}$$

Next we need `functional controls': this is a   {\it  cocone with
base $D$ and vertex an object $C$}.
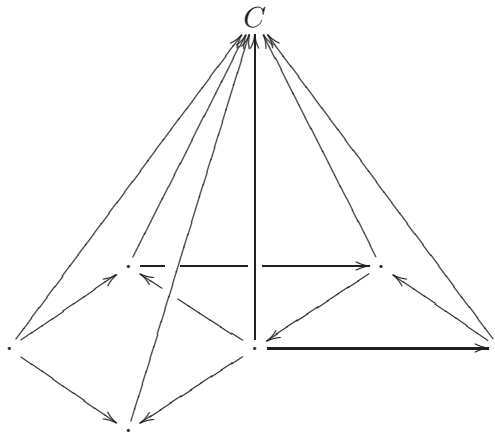
$$\xymatrixcolsep{3pc} \xymatrixrowsep{2pc}\xymatrix{
&&C&&\\&&&&\\
&&&&\\
 & . \ar@{-}[r]|(0.35)\hole\ar[uuur]&\ar[r]&.\ar[dl]\ar[uuul]&\\
.\ar[ur]\ar[uuuurr]\ar[dr]& & .
\ar[ul]|!{[dl];[uuuu]}\hole\ar[rr]\ar[dl]\ar[uuuu]
 & &.\ar[ul]\ar[uuuull]\\
&.\ar[uuuuur]&&&}$$ such that  each of the triangular faces of
this cocone is commutative.

 The output from such input data will be
an object $\mathsf{colim}( D)$ in our category $\C$ defined by a
special {\em colimit
 cocone}  such that any cocone on $D$ factors {\it uniquely} through the colimit
 cocone. The commutativity condition on the cocone in essence
 forces,
 in the colimit, an interaction of the images of different parts of the diagram
 $D$. The uniqueness condition makes the colimit the {\it best
 possible} solution to this factorisation problem

In the next picture  the colimit is written $$\blacklozenge=
\mathsf{colim}(D),$$the dotted arrows represent new morphisms
which combine to make the colimit cocone:
$$\xymatrixcolsep{3pc} \xymatrixrowsep{2pc}
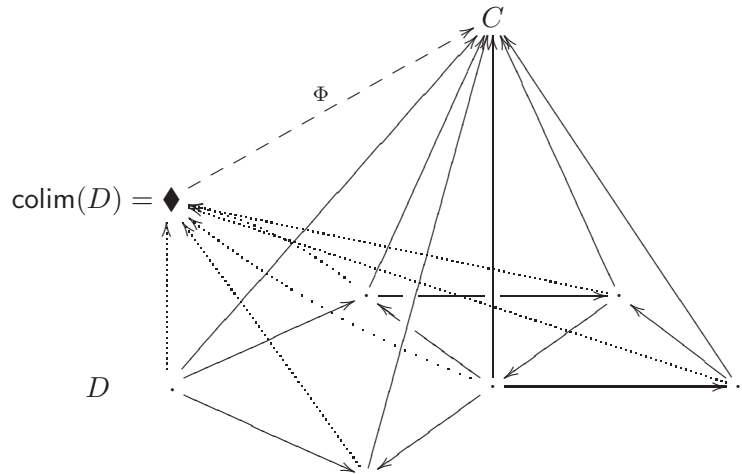
\xymatrix{
&&C&&\\&&&&\\
**[l]\mathsf{colim}(D)=\blacklozenge\ar@{-->}[rruu]^ \Phi &&&&\\
 & .\ar@/_/@{.>}[lu] \ar@{-}[r]|(0.35)\hole\ar[uuur]&\ar[r]&.\ar@{.>}[lllu]\ar[dl]\ar[uuul]&\\
**[l]D \qquad \cdot  \ar@{.>}[uu] \ar[ur]\ar[uuuurr]\ar[dr]& & .
\ar@/^/@{.>}[lluu] \ar[ul]|!{[dl];[uuuu]}\hole
 \ar[rr]\ar[dl]\ar[uuuu]
 & &.\ar[ul]\ar@{.>}[uullll]\ar [uuuull]\\
&.\ar[uuuuur] \ar@{.>}[uuul]&&&}$$  and the broken arrow $\Phi$ is
constructed from the other information. Again, all triangular
faces of the combined picture are commutative. Now stripping away
the `old' cocone gives the factorisation of the cocone via the
colimit:
$$\xymatrixcolsep{3pc} \xymatrixrowsep{2pc}
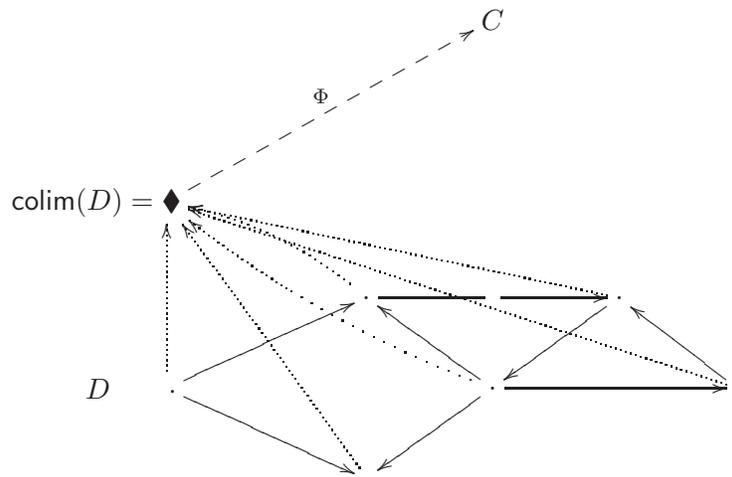
\xymatrix{
&&C&&\\&&&&\\
**[l]\mathsf{colim}(D)=\blacklozenge  \textcolor{red}{\ar@{-->}[rruu]^\Phi} &&&&\\
 & .\ar@/_/@{.>}[lu] \ar@{-}[r]&\ar[r]&.\ar@{.>}[lllu]\ar[dl]&\\
**[l]D \qquad \cdot  \ar@{.>}[uu] \ar[ur]\ar[dr]& & .
\ar@/^/@{.>}[lluu] \ar[ul]\ar[rr]\ar[dl]
 & &.\ar[ul]\ar@{.>}[uullll]\\
&. \ar@{.>}[uuul]&&&}$$

\noindent {\bf Intuitions:}

The object colim($D$) is `put together' from, or `composed of the
parts of',  the constituent diagram $D$ by means of the colimit
cocone. From beyond (or above
 our diagrams) $D$, an object $C$ `sees' the diagram $D$
`mediated' through its colimit, i.e. if $C$ tries to interact with
the whole of $D$, it has to do so via colim($D$). The colimit
cocone can be thought of as a kind of program: given any cocone on
$D$ with vertex $C$, the output will be a morphism $$\Phi:
\mathsf{colim}(D)\to C$$ constructed from the other data. How is
this morphism realised, what are its values?

To focus on a common example, consider the process of sending an
email document, call it $E$. To send this we need a server $S$,
which breaks down the document $E$ into many parts $E_i$ for $i$
in some indexing set $I$, and labels each part $E_i$ so that it
becomes $E_i'$. The labelled parts $E'_i$ are then sent to various
servers $S_i$ which then send these as messages $E''_i$ to a
server $S_C$ for the receiver  $C$. The server $S_C$ combines the
$E''_i$ to produce the received message $M_E$ at $C$. Notice also
that there is an arbitrariness in breaking the message down, and
in how to route through the servers $S_i$, but the system is
designed so that the received message $M_E$ is independent of all
the choices that have been made at each stage of the process. A
description of the email system as a colimit may be difficult to
realise precisely, but this analogy does suggest the emphasis on
the amalgamation of many individual parts to give a working whole,
which yields exact final output from initial input, despite
choices at intermediate stages.

One question for neuroscientists is: does  the brain use analogous
processes for communication between its various structures? What
we can say is that this general colimit notion represents a
general mathematical process which is of fundamental importance in
describing and calculating with many algebraic and other
structures.

To go back to our email analogy,  the morphism $\Phi$ is
constructed on some element $z$ of  colim($D$), by splitting $z$
up somehow into pieces $z_i$ which come from parts $z''_i$ in
objects of $D$, mapping these $z''_i$ over using the cocone on
$D$, and combining them in $C$. This is how some proofs that
certain cocones are colimits are actually carried out, see
\cite{Brown:fields}.

Thus a conjecture as far as biological processes are concerned, is
that this notion of colimit may give useful analogies to the way
complex systems operate. More generally, it seems possible that
this particular concept in category theory, seeing how a big
object is built up of smaller related pieces, may be useful for
the mathematics of processes.

\section{Higher dimensional algebra}
The aim here is to explain some mathematical ideas with which the
authors have been preoccupied since the 1960s and 1970s
respectively. There was a lot of experimentation to produce the
mathematics which would encompass some apparently simple
intuitions. This experimentation can be viewed as `extraction of
concepts', and it also exemplifies mathematical concepts that
might provide a model of inter-relationships that is much freer
than the usual ones. Indeed, higher dimensional algebra is being
used to model distributed systems.

A basic idea is that we may need to get away from `linear'
thinking in order to express intuitions clearly.

Thus the equation \begin{equation}2 \times (5+3)=2 \times 5 + 2
\times 3 \label{formula}\end{equation}  is more clearly shown by
the figure \begin{equation} \label{picture} \begin{split}|\;|\;|\;|\;|&\qquad |\;|\;|\\
|\;|\;|\;|\;|&  \qquad |\;|\;| \end{split} \end{equation} Indeed
the number of conventions you need to understand equation
\eqref{formula} make it seem barbaric compared with the picture
\eqref{picture}. It is also interesting to see how you could
express in a picture the general linear formula of the
distributivity law
$$a \times (b+c)= a \times b + a \times c.$$  The
importance of having simple comprehensible pictures instead of
complex formulae is that the pictures help one to imagine theorems
and their proofs. (The above exposition is borrowed from an
account in week 53 of John Baez's series `This week's find in
mathematical physics' \cite{Baez}.)

Those who have read Edwin Abbott's famous book `Flatland'
\cite{abbott} (and those who have not, have a delight in store!)
will recall the limited interactions available to the inhabitants
of Lineland. It seems unreasonable to suppose that a purely linear
mathematics can express reasonably the complex interactions that
occur in the brain.

We often translate geometry into algebra. For example, a figure as
follows:
$$
\xymatrix@M=0pc { \bullet \ar @{-}[r] |@{>} ^a & \bullet \ar
@{-}[r] |@{>}^b & \bullet
 \ar @{-}[r] |@{>}^c & \bullet \ar @{-}[r] |@{>} ^d & \bullet}
$$
is easily translated into
$$abcd$$
and the language for expressing this is again that of {\em
category theory}. It is useful to express this intuition as
`composition is an algebraic inverse to subdivision'. The labelled
subdivided line gives the composite word, $abcd$.

But how do we express a diagram such as the following
$$ {\objectmargin{0.1pc} 
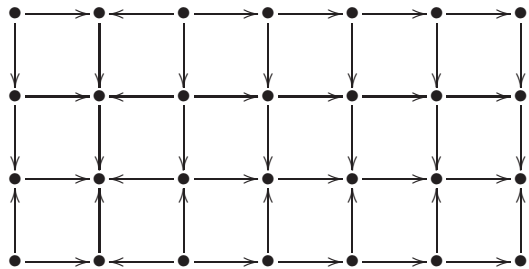
\diagram \bu \rto \dto &
\bu  \dto & \lto \bu \rto \dto & \bu
\rto \dto & \bu \rto \dto & \bu \rto \dto & \bu  \dto \\
 \bu\rto \dto &  \bu  \dto & \lto\bu \rto \dto & \bu \rto \dto & \bu
\rto \dto & \bu \rto \dto & \bu  \dto \\
\bu \rto  &  \bu & \lto\bu \rto  &  \bu \rto  & \bu \rto  & \bu \rto & \bu \\
\bu \rto \uto  & \bu \uto   &\lto  \bu \rto \uto   & \bu \rto
\uto
 & \bu \rto   \uto & \bu \rto \uto  & \bu  \uto
\enddiagram }$$
where the squares are supposed filled and labelled? It seems that
in doubling the number of dimensions from 1 to 2,   you need to
move from categories to {\em double categories}, or something
similar, based on directed squares rather than on arrows.

The extra richness involved is that a square can have more
complicated relations to other squares than can happen in the
linear situation.

Such questions arose from a {\em gluing} or colimit problem in
topology, namely to describe the behaviour of a big object in
terms of the behaviour of its parts.  When the first author
started work on this, the particular problem in dimension 1 was
well known and solved, but the question was to carry out similar
methods in higher dimensions. For a survey of this for a
mathematical audience, see \cite{Brown:fields}.

In the  topological situation from which these ideas arose, to
obtain the uniqueness in the output of a colimit, as we always
require for our emails,  we had to go further than {\it `algebraic
inverses to subdivision'} and to use also the notion of {\it
`commutative cube'}.

A commutative square is very easy:

 $$\def\labelstyle{\textstyle}  
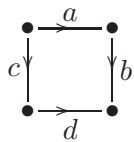
\xymatrix@M=0.1pc @=2pc{\bullet \ar @{-} [r]|@{>} ^a  \ar @{-} [d]|@{>} _c & \bullet \ar @{-}
 [d]|@{>} ^b \\ \bu  \ar @{-} [r]|@{>} _d & \bu}
$$  is easily translated into mathematics as \begin{equation} ab=cd, \qquad \mbox{or} \qquad a= cdb^{-1}.
\end{equation}

The surprising thing is that to determine a {\it commutative cube}
needs some new ideas. First we need to know how to compose the
square faces of a cube. This can be done in two directions:

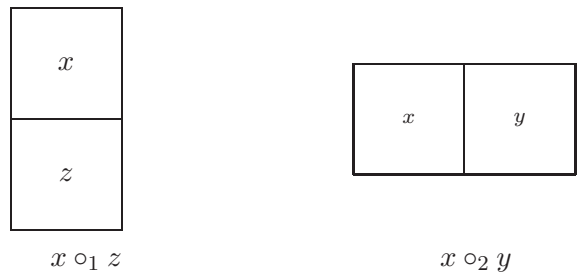
\begin{align*} { \def\labelstyle{\textstyle} \objectmargin{0pc} \spreaddiagramrows{1.5pc}
\spreaddiagramcolumns{1.5pc}\diagram \midsq{x} \rline \dline &
\dline
\\ \midsq{z} \rline \dline & \dline \\ \rline &
\enddiagram } \qquad\qquad &\qquad\qquad{\objectmargin{0pc}\spreaddiagramrows{1.5pc}
\spreaddiagramcolumns{1.5pc} \diagram \midsq{x} \rline \dline &
\midsq{y}  \rline \dline  &  \dline \\ \rline & \rline &
\enddiagram   } \\ \hspace{-2em} x\circ_1 z \qquad \qquad & \hspace{7em}   x
\circ_2 y
\end{align*}
So we get a notion of a {\bf double category}, whose elements are
squares, for which in the above diagram the composition $x
\circ_1z$ is defined if and only if the bottom edge of $x$ is the
same as the top edge of $z$, and similarly (but right and left
edges) for $x \circ_2 y$.  Thus the  compositions are partially
defined, under geometric conditions. This enables a close relation
between the algebra and geometry. On these compositions, we have
to impose all the obvious geometric rules.

These rules  enables an easy description of multiple compositions.
But there remains the question of how to  define a {\it
commutative cube}? A cube has six faces, which can divide into two
groups of three. With clear conventions, one would like to equate
the two compositions
\begin{equation}\def\labelstyle{\textstyle} {{\objectmargin{0pc}\spreaddiagramrows{1.5pc}
\spreaddiagramcolumns{1.5pc}
 \diagram \midsq{\prt^0_2} \rline \dline & \midsq{\prt^1_3}
 \rline \dline & \dline \\
\midsq{\prt^1_1} \rline \dline &  \rline \dline & \\
\rline &  & \enddiagram } \hspace{6em}
{\objectmargin{0pc}\spreaddiagramrows{1.5pc}
\spreaddiagramcolumns{1.5pc}
\diagram   & \midsq{\prt^0_1} \rline \dline & \dline \\
\midsq{\prt^0_3} \rline \dline & \midsq{\prt^1_2} \rline \dline & \dline \\
\rline & \rline & \enddiagram} } \end{equation}Unfortunately, this
does not make sense because the little squares do not form a
rectangular array, and also because the edges of each block are
not correct to give equality. It is interesting that the extension
from squares in dimension 2 to cubes in dimension 3 produces these
gaps which require a new set of concepts to handle them.

We need new elements to fill in the corners, and in fact you also
need to expand out, to obtain the equality of:
\begin{equation}\def\labelstyle{\textstyle}
{{\objectmargin{0pc}\spreaddiagramrows{1.5pc}
\spreaddiagramcolumns{1.5pc}\diagram \midsq{\hh} \rline \dline &
\midsq{\prt^0_2} \rline \dline &\midsq{\prt^1_3} \dline \rline& \dline \\
\midsq{\tl} \rline \dline & \midsq{\prt^1_1} \rline \dline
&\midsq{\br}\rline \dline
&\dline \\
\rline & \rline & \rline &  \enddiagram}   \hspace{6em}
{\objectmargin{0pc}\spreaddiagramrows{1.5pc}
\spreaddiagramcolumns{1.5pc}
 \diagram \midsq{\tl} \rline \dline &
\midsq{\prt^0_1} \rline \dline &\midsq{\br} \dline \rline& \dline \\
\midsq{\prt^0_3} \rline \dline & \midsq{\prt^1_2} \rline \dline
&\midsq{\hh}\rline \dline
&\dline \\
\rline & \rline & \rline &  \enddiagram} }\end{equation} Thus
2-dimensional algebra needs some new basic constructions. This is
not surprising. In dimension 1, you are limited to staying still,
moving forward, or moving backward. In dimension 2, you can also
turn left or right. This is what needs to be modelled formally.

We need  some  special squares called {\em thin squares}:
\begin{center}
\begin{tabular}{ccc}
$\quads{1}{1}{1}{1}$& $\quads{a}{1}{a}{1}$& $\quads{1}{b}{1}{b}$
\\
&&\\
$\tsq$ & $\hh \mbox{ or }\eps_2 a$ & $\vv \mbox{ or }\eps_1 b$
\end{tabular}
\end{center}
of which the above first three are forms of horizontal and
vertical identities, which correspond to not moving in certain
directions, while our last two
\begin{center}
\begin{tabular}{cc}
$\quads{a}{a}{1}{1}$& $\quads{1}{1}{a}{a}$  \\&\\
 $\br \mbox{ or }\Gamma a$ & $\tl\mbox{ or } \Gamma' a$
\end{tabular}
\end{center}
are known as {\em connections}, and correspond to turning left or
right. How do these all interact?

 The rules on the connections are as follows:
\begin{enumerate}[\rm (i)]
\item $ [\, \tl \;\br \;] = \vv$,\quad $[\,\Gamma'a \; \Gamma a
\;]= \eps_1 a$; \item $ \begin{bmatrix}\,\tl\,\\
\br\,\end{bmatrix}  = \hh\;$, \quad $ \begin{bmatrix}\,\Gamma' a\,\\
\Gamma a\,\end{bmatrix} = \eps_2 a$; \item $
\begin{bmatrix}\tl & \hh \\ \vv & \tl
\end{bmatrix} = \tl $,\quad$ \begin{bmatrix}\Gamma'a  & \eps_2 a \\
\eps_1 a & \Gamma' b
  \end{bmatrix} = \Gamma' (ab)$;\quad
\item $ \begin{bmatrix}\br & \vv \\ \hh & \br  \end{bmatrix} = \br
\,$, \quad$ \begin{bmatrix}\Gamma a  & \eps_1 b  \\ \eps_2 b &
\Gamma b
\end{bmatrix} = \Gamma (ab)$.
\end{enumerate}
Notice that we give two notations for these thin squares, and that
the more heiroglyphic notation is much easier to grasp by the eye.
This emphasises the idea that the history of mathematics is much
involved with the history of improved notation. The first two
rules for connections can be thought of as formalising  turning
left and then right (or the other way round), while the last two
formalise that turning with an arm outstretched is the same as
just turning.

Here we just have time to show a typical calculation to convince
you that 2-dimensional rewriting  looks like a new kind of
manipulation. We start with the first complex diagram and use the
rules to rewrite it to a simpler one. The rewriting may also be
done a different way to give a different expression which has, by
the rules, an equal evaluation in the double category with
connections.
$$\def\labelstyle{\textstyle} \spreaddiagramrows{1pc} \spreaddiagramcolumns{1pc}
\objectwidth{0in} \objectmargin{0in} \def\labelstyle{\textstyle}
\diagram \rline \dline \midsq{\sq} & \rline \dline \midsq{\sq}
   & \rline \dline \midsq{\sq} &
\rline^c \dline \midsq{\vv} & \rline \dline \midsq{\tl} & \dline^d \\
\rline \dline \midsq{\tl} & \rline \dline \midsq{\hh} & \rline
\dline^a &
\rline \dline_c \midsq{\br} & \rline \dline \midsq{\vv} & \dline \\
\rline \dline \midsq{\vv} & \rline \dline \midsq{\tl} & \dline^b
\save \go+<1.5pc,0pc> \Drop{\hh} \restore &
\rline \dline_d \midsq{\hh} & \rline \dline \midsq{\br} & \dline \\
\rline \dline_a \midsq{\br} & \rline \dline \midsq{\vv}
   & \rline \dline \midsq{\sq} &
\rline \dline \midsq{\sq} & \rline \dline \midsq{\sq} & \dline \\
\rline & \rline_b & \rline & \rline & \rline &
\enddiagram
$$
$$
= \def\labelstyle{\textstyle} \spreaddiagramrows{1pc}
\spreaddiagramcolumns{1pc} \objectwidth{0in} \objectmargin{0in}
\def\labelstyle{\textstyle} \diagram \rline \dline & \rline &
\rline \dline \midsq{\sq} &
\rline^c \dline \midsq{\vv} & \rline \dline \midsq{\tl} & \dline^d \\
\rline \dline & \rline \save \go+<0pc,1.5pc> \Drop{\sq}
   \restore & \rline \dline_{ab} \midsq{\hh} &
\rline_c \dline^{cd} & \rline_d & \dline \\
\rline^a \dline_a \midsq{\br} & \rline^b \dline \midsq{\vv} \save
\go+<0pc,1.5pc> \Drop{\tl} \restore & \rline \dline \midsq{\sq} &
\rline \dline & \rline \save \go+<0pc,1.5pc> \Drop{\br} \restore & \dline \\
\rline & \rline_b & \rline & \rline & \rline \save \go+<0pc,1.5pc>
   \Drop{\sq} \restore &
\enddiagram
= \diagram
\rline^c \dline \midsq{\vv} & \rline \dline \midsq{\tl} & \dline^d \\
\rline_c \dline & \rline_d & \dline \\
\rline^a \dline_a \midsq{\br} & \rline^b \dline \midsq{\vv} \save
\go+<0pc,1.5pc> \Drop{\vv} \restore
   & \dline \\
\rline & \rline_b &
\enddiagram   $$
Analogous rewriting has been carried out in three dimensions in
\cite{aa-b-s}, and, as may be imagined, is not easy to handle. As
always, the development of a new mathematics solves some problems
and then brings a range of new problems into view.

Rewriting, changing one formula to another according to certain
rules,  is a basic facet  of much mathematics. When multiple
rewrites are occurring,  for various operations,  the study of
this `distributed' or `concurrent' rewriting  can often be
pictured as involving higher dimensional laws. Models of the
complex message passing and distributed rewriting in the brain
have yet to find the laws relevant to that context. Such laws
could be key elements for the next level of understanding of
cognition.

{{\bf Questions}}

Is the colimit notion useful to describe the way the brain (or a
brain module)  integrates structural information incoming in
various forms to give a determined output?

What other models are there? It is good to start with the
assumption that the simplest idea works!

Information is often `subdivided' by the sensory organs and is
reintegrated by the brain. To  enable different parts of that
information to be integrated, there must be some `glue', some
inter-relational information available. If we are given arrows
$a,b,c,d$ with no information on where they start or end, then we
could form combinations which make no geometric sense. The
colimit/composition process makes sense only where the
inter-relations are also such as to enable the `integration' to be
well defined. Higher dimensional algebra allows more complex
notions of `well formed composition', and ones more adapted to
geometry.

{\bf Computation and computer science}

A mathematician always wants to do sums, get explicit answers to
some situations, but recognises that not every sum can be worked
out. In fact you may need more mathematics to work out how to do
the sums.

Computer languages do not yet seem to be good at doing {\it
structural mathematics}, or even at expressing a high level
algorithm, i.e. one which works at a high structural level.
Indeed, current computer languages were not designed with the
needs of mathematics in mind. Yet  algorithms designed to encode a
high level of structural information should be more efficient --
we do not tell people the way to the station by giving information
on  all the cracks in the pavement. We need to map the landscape,
that is give a usable model of it,  to describe the path through
it.

Since evolution is `concerned with' efficiency, we must expect
that the brain has evolved methods for dealing with structural
information. It looks like a reasonable conjecture that category
theory and higher dimensional category theory could be necessary
for modelling this kind of behaviour.

Higher Dimensional Algebra has already shown its use in models of
information management, and in concurrency. Descriptions of
systems by graphs are well known, with development described
algebraically by paths in graphs.  Interacting systems need higher
dimensional graphs, and a generalisation of the notion of path.
Higher Dimensional Algebra is still young, and there are many new
possibilities opened up, as a web search shows.

\end{document}